\documentclass[10pt]{article}
\usepackage{amssymb,amsmath,textcomp}
\newtheorem{satz}{Theorem}[section]
\newtheorem{lemma}{Lemma}[section]

\newtheorem{bemerk1}{Remark}[section]
\newtheorem{bei}{Example}[section]
\newtheorem{korollar}{Corollary}[section]

\newcommand{\iR}{\mathbb{R}}
\newcommand{\iN}{\mathbb{N}}

\newcommand{\iC}{\mathbb{C}}

\newcommand{\oH}{\hspace*{0.39em}\raisebox{0.6ex}{\textdegree}\hspace{-0.72em}H}
\DeclareMathOperator*{\esup}{ess\,sup}
\DeclareMathOperator*{\einf}{ess\,inf}

\begin{document}
\begin{center}
{\bf\Large Stability, instability, and blowup for time fractional and other non-local in time semilinear subdiffusion equations}
\end{center}
\vspace{0.5em}
\begin{center}
Vicente Vergara\footnote{V.V.\ was partially supported by FONDECYT grant 1150230.} and Rico Zacher\footnote{R.Z.\ was partially supported by a Heisenberg fellowship of the German Research Foundation (DFG), GZ  Za 547/3-1.}
\end{center}
\vspace{0.3em}
\begin{center}
{\em Dedicated to Jan Pr\"uss on the occasion of his 65th birthday}
\end{center}

\begin{abstract}
We consider non-local in time semilinear subdiffusion equations on a bounded domain, where the kernel in the integro-differential operator
belongs to a large class, which covers many relevant cases from physics applications, in particular the important case of fractional
dynamics. The elliptic operator in the equation is given in divergence form with bounded measurable coefficients.
We prove a well-posedness result in the setting of bounded weak solutions and study the stability and instability of the zero
function in the special case where the nonlinearity vanishes at $0$. We also establish a blowup result for positive convex and
superlinear nonlinearities.   
\end{abstract}
\vspace{0.7em}
\begin{center}
\end{center}

\section{Introduction}
Let $\Omega$ be a bounded domain in $\iR^N$. We consider the problem
\begin{align}
\partial_t\big(k\ast[u-u_0]\big)-\mbox{div}\big(A(t,x)\nabla u\big) & =f(u),\quad t>0,\,x\in\Omega,
\nonumber\\
u & =0,\quad t>0,\,x\in \partial\Omega, \label{subproblem}\\
u|_{t=0} & = u_0,\quad x\in \Omega.\nonumber
\end{align}
The kernel $k\in L_{1,\,loc}(\iR_+)$ is given, and $k\ast v$ denotes the
convolution on the positive halfline $\iR_+:=[0,\infty)$ w.r.t.\ the time variable,
that is
$(k\ast v)(t)=\int_0^t k(t-\tau)v(\tau)\,d\tau$, $t\ge 0$. Note that for sufficiently smooth $u$ with $u(0)=u_0$, 
\begin{equation} \label{operatorreform}
\partial_t \big(k\ast [u-u_0]\big)=k \ast \partial_t u.
\end{equation}

The kernel $k$ belongs to a large class of kernels, it is merely assumed to satisfy the condition
\begin{itemize}
\item [{\bf ($\mathcal{PC}$)}] $k\in L_{1,\,loc}(\iR_+)$ is nonnegative and nonincreasing, and there exists a kernel $l\in L_{1,\,loc}(\iR_+)$ such that
$k\ast l=1$ on $(0,\infty)$.
\end{itemize}
In this case we say that $k$ is a kernel of type $\mathcal{PC}$ (cf. \cite{ZWH}) and also write $(k,l)\in \mathcal{PC}$. Note that $(k,l)\in {\cal PC}$ implies that $l$ is completely
positive, cf.\ \cite[Theorem 2.2]{CN} and \cite{CN1}, in particular $l$ is nonnegative.

Condition ($\mathcal{PC}$) covers most of the relevant integro-differential operators w.r.t.\ time that
appear in physics applications in the context of subdiffusion processes. An important example is given by $(k,l)=(g_{1-\alpha},g_\alpha)$ with $\alpha\in(0,1)$, where
$g_\beta$ denotes the standard kernel
\[
g_\beta(t)=\,\frac{t^{\beta-1}}{\Gamma(\beta)}\,,\quad
t>0,\quad\beta>0.
\]
In this case, the term $\partial_t(k\ast v)$ becomes the classical Riemann-Liouville fractional derivative
$\partial_t^\alpha v$ of order $\alpha$, and $k \ast \partial_t v={}^c D_t^\alpha v$, the Caputo fractional derivative (cf.\
the right-hand side in (\ref{operatorreform})), of the (sufficiently smooth) function $v$, see e.g.\ \cite{KST}.

Another interesting example is given by the pair
\begin{equation} \label{ultrapair}
k(t)=\int_0^1 g_\beta(t)\,d\beta,\quad l(t)=\int_0^\infty \,\frac{e^{-st}}{1+s}\,{ds},\quad t>0.
\end{equation}
In this case the operator $\partial_t(k\ast \cdot)$ is a so-called operator of
{\em distributed order}, see e.g.\ \cite{Koch08,VZ}. Further examples will be discussed in Example \ref{kerne} below.

Concerning the
coefficients $A=(a_{ij})$ we assume
that
\begin{itemize}
\item [{\bf ($\mathcal{H}$)}] $A\in L_\infty((0,T)\times \Omega;\iR^{N\times
N})$ for all $T>0$, and $\exists \nu>0$ such that
\[
\big(A(t,x)\xi|\xi\big)\ge \nu|\xi|^2,\quad \mbox{for
a.a.}\;(t,x)\in (0,\infty)\times \Omega,\,\mbox{and all}\,\xi\in
\iR^N.
\]
\end{itemize}

Problems of the form (\ref{subproblem}) with $f=0$, in particular time fractional diffusion equations,
have attracted much interest during the last years, mostly due to their applications in the modeling of anomalous diffusion, see
e.g.\ \cite{Koch08, Koch11, Metz, Uch} and the references therein for the physical background.
To provide some more specific motivation, consider for the moment the case $\Omega=\iR^N$, $A(t,x)=I$ and $f=0$, and
let $Z(t,x)$ denote the fundamental solution of the corresponding equation satisfying  $Z|_{t=0}=\delta_0$. 
If $(k,l)\in \mathcal{PC}$, then $Z$ can be constructed via subordination from the heat kernel and one can show that $Z(t,\cdot)$ is a probability density function on $\iR^N$ for all $t>0$, see \cite{KSVZ}. Further, 
the so-called {\em mean
square displacement}, which is defined in our case as
\[
m(t)=\int_{\iR^N}|x|^2 Z(t,x)\,dx,\quad t>0,
\]
and describes how fast particles diffuse, is known to be given by
\[
m(t)= \,2N\,(1\ast l)(t),\quad t> 0,
\]
see \cite{KSVZ}. In the time fractional diffusion case (i.e.\ the first
example) one observes that $m(t)=ct^\alpha$ with some constant $c>0$ (see also \cite{Metz}), which shows that the diffusion is slower than in the classical case of Brownian motion, where $m(t)=c t$. In our second example, the mean square displacement $m(t)$ behaves like $c\log t$ for $t\to \infty$, see \cite{Koch08}. In this case
the corresponding diffusion equation describes a so-called {\em ultraslow diffusion} process.

Semilinear problems of the form (\ref{subproblem}) generalize the pure diffusion case (with $f=0$) by including a
nonlinear source term. Such problems also occur as models for nonlinear heat flow in materials with memory, see e.g.\
\cite{GP, Nun}. 

The main results of this paper are the following. Assuming that $f$ is locally Lipschitz continuous we first establish local well-posedness for (\ref{subproblem})
in the framework of bounded weak solutions. We also show that there is a maximal interval of existence $[0,t_*)$ and 
that the solution blows up as $t\to t_*-$ if $t_*<\infty$ and $f$ is defined on all of $\iR$, see Theorem 
\ref{wellposedsubdiffusion} below. We point out that in our local
well-posedness result the coefficient matrix $A$ is only assumed to satisfy condition ($\mathcal{H}$), in particular it is allowed to depend on time. 

The second main result, Theorem \ref{stabilitysubdiffusion}, 
provides sufficient conditions for the stability and asymptotic stability of
the zero function; here we assume that $f(0)=0$ and that $f'(0)$ exists. In addition we restrict ourselves to the
case where $A$ is independent of $t$. In the special case, where $A$ is also symmetric, the stability condition is
given by $f'(0)<\lambda_*$, where $\lambda_*>0$ denotes the smallest eigenvalue of the operator
$\mathcal{L}v=-\mbox{div}\big(A(x)\nabla v\big)$ (with Dirichlet boundary condition) in $L_2(\Omega)$. 
In the symmetric case we also prove instability of the zero function if $f'(0)>\lambda_*$ and $l\notin L_1(\iR_+)$,
see Theorem \ref{instabilitysubdiffusion}.

The last main result, Theorem \ref{blowupPDE}, is concerned with the blowup of solutions to (\ref{subproblem}).
Here we assume that $A$ is independent of $t$ and symmetric. Concerning $f$ we impose a convexity and superlinearity
condition; these conditions on $f$ are also widely used in the classical theory of parabolic PDEs. We prove that for
sufficiently large initial data $u_0\ge 0$ the corresponding weak solution to  (\ref{subproblem}) blows up in finite time.
We remark that in contrast to the instability result, here the kernel $l$ is also allowed to be integrable on $\iR_+$.

Our proofs of the well-posedness and stability results require a couple of auxiliary results such as, e.g., a comparison
principle and an appropriate linear stability result. Some of these results seem to be new and are interesting in its own
right.    
  
In view of condition ($\mathcal{PC}$) the nonlocal PDE in (\ref{subproblem}) can be rewritten as a Volterra
equation on the positive halfline with a completely positive kernel; this can be seen by convolving the PDE with the kernel
$l$. If $A$ is independent of $t$, the problem can be viewed as an abstract Volterra equation of the form
\begin{equation} \label{Voltab}
v(t)+(l\ast \mathcal{L}v)(t)=v_0+(l\ast F(v))(t),\quad t\ge 0.
\end{equation}
Here $v$ takes values in some Banach space of functions of the spatial variable and $\mathcal{L}$ denotes the elliptic
operator mentioned above. There has been a substantial amount of work on such abstract (linear and nonlinear) Volterra and integro-differential equations since the 1970s,
in particular on existence and uniqueness, regularity, and long-time behaviour of solutions,
see, for instance, \cite{AP, CLS, CNa, CN, Grip2,Grip1, ZEQ}, and the monograph \cite{JanI}. The results in \cite{CLS} contain as
a special case the local well-posedness in continuous interpolation spaces of an abstract time fractional quasilinear equation
and are applicable to (\ref{Voltab}) with $l=g_\alpha$, $\alpha\in (0,1)$. The existence and uniqueness results in \cite{Grip1}
are based on the theory of accretive operators and also apply to (\ref{Voltab}). 

However, these abstract results are not
applicable to solve (\ref{subproblem}) with rough coefficient matrix $A(t,x)$. In order to achieve this, an appropriate
theory of weak solutions is required. In this paper we make use of the results from \cite{ZWH} on weak solutions
to abstract evolutionary
integro-differential equations in a Hilbert space setting. In the case with rough coefficient matrix $A(t,x)$ and $f=0$,
optimal $L_2$-decay estimates were proved in \cite{VZ}. 

In \cite{CP}, the authors establish
the global existence in a strong $L_p$-setting for a semilinear parabolic Volterra equation with Dirichlet boundary condition; here $A(t,x)=I$, but the nonlinearity
may also depend on $t,x$ and $\nabla u$.     
        
As to stability results in the nonlinear case, linearized stability for an abstract Volterra equation has been studied in \cite{Kato} in the abstract framework
of accretive operators. Due to the assumptions on the kernel, the results in \cite{Kato} do not apply to our situation.
Even the important time fractional case is excluded there.

Concerning blowup results, there exist already some results on special cases of  (\ref{subproblem}), mostly with
fractional dynamics. In \cite{AAAKT}, the authors show blowup in the case $f(x)=x^2-x$ and $k=g_{1-\alpha}$,
thereby solving a problem which was raised in \cite{NSY}.
Time fractional
semilinear diffusion equations with power type nonlinearity in the whole space are studied in \cite{KLT}. The authors
determine critical exponents of Fujita type and establish necessary and sufficient conditions for the existence of nontrivial
global solutions. Blowup for one-dimensional time fractional diffusion problems with a source term of the 
special form $g_{1-\alpha}\ast (\delta(x-a)f(u(\cdot,a))(t)$, where $a\in \Omega$ is fixed and $\delta$ denotes the
Dirac delta distribution is studied in \cite{OR}. In the purely time-dependent case, that is, in the case without elliptic operator, one finds many blowup results 
for several kinds of Volterra equations, see e.g.\  \cite{MO,RLO}. However the general situation we consider in Theorem
\ref{basicblowup} with a kernel $k$ of type $\mathcal{PC}$ and a general positive, convex and superlinear
function $f$ does not seem to have been studied so far. Moreover, the argument we give to prove Theorem \ref{basicblowup},
although well known in the ODE case,
seems to be totally new in the context of Volterra equations.

We would like to point out that the theory developped in this paper can be extended to a more general class
of equations which is obtained by adding an additional term $k_0 \partial_t u$, with $k_0\ge 0$, on the left-hand side of (\ref{subproblem}). In this situation, one has to replace condition ($\mathcal{PC}$) by
\begin{itemize}
\item [{\bf ($\mathcal{PC}'$)}]
$k\in L_{1,\,loc}(\iR_+)$ is nonnegative and nonincreasing, and there exists a kernel $l\in L_{1,\,loc}(\iR_+)$ such that
$k_0 l(t)+(k\ast l)(t)=1$ on $(0,\infty)$.
\end{itemize}
Also in this situation, $l$ is completely
positive (see e.g.\ \cite{CN}), and thus nonnegative. The special choice $k_0=1$, $k=0$ and $l=1$ in this more general formulation leads to the classical parabolic case, which is not covered by (\ref{subproblem}) under assumption ($\mathcal{PC}$). We believe that the analogue (in the classical parabolic case) of our nonlinear well-posednes result, Theorem \ref{wellposedsubdiffusion}, is known,
however we could not find a reference. It is not difficult to check that our arguments can be generalized
to the case with additional term $k_0 \partial_t u$, since $l$ is the important kernel. However, this also requires an extension of some auxiliary results cited from the literature in Section 3, in particular of the basic linear existence result, Theorem \ref{linearweakex}. This is possible but not the subject of this paper,
so we confine ourselves to the situation described above.

As to other boundary conditions, we are convinced that our arguments can be adapted to obtain corresponding results for
homogenous Neumann and Robin boundary conditions. Again, to achieve this one also has to prove first analogues of some
of the auxiliary results cited from the literature in Section 3. 

Concerning limitations of our theory, we remark that in our framework, the nonlinearity $f$ cannot be allowed to depend on the gradient $\nabla u$. Our method does not seem to apply either to quasilinear equations, that is to the situation where the
coefficient matrix $A$ in (\ref{subproblem}) also depends on $u$ or $\nabla u$. This is due to the lack of regularity.
Note, however, that for $A=A(u)$ and $k=g_{1-\alpha}$ with $\alpha\in (0,1)$ and assuming more regularity on the
initial value $u_0$, the corresponding quasilinear problem can be solved uniquely by means of maximal regularity, see
e.g.\ \cite{CLS,ZaG}.   

Another limitation is that superdiffusion
equations are excluded, in particular the time fractional case with time order $\alpha\in (1,2)$, by the lack
of the maximum principle. Equations of this type have also been studied quite intensively, we refer to 
\cite{CLS, Fu, SW, JanI, PVZ, ZEQ, ZQL} and the references given therein.   

The paper is organized as follows. In Section 2 we collect basic properties of kernels of type $\mathcal{PC}$ and state
a fundamental convexity property for operators of the form $\frac{d}{dt}(k\ast \cdot)$. Section 3 is devoted to 
linear problems and provides important tools, which are required in the analysis of the nonlinear problem.
Section 4 contains the well-posedness result for (\ref{subproblem}). Section 5 deals with stability and instability of
the zero function. Finally, in Section 6 we study the blowup of solutions, first in the purely time-dependent case and
then in the full PDE case.

\section{Preliminaries} \label{Prelim}
We first collect some properties of kernels of type $\mathcal{PC}$. Let $(k,l)\in \mathcal{PC}$.
For $\gamma\in \iR$ define the kernels $s_\gamma, r_\gamma \in L_{1,loc}(\iR_+)$ via the scalar Volterra equations
\begin{align*}
s_\gamma(t)+\gamma(l\ast s_\gamma)(t) & = 1,\quad t>0,\\
r_\gamma(t)+\gamma(l\ast r_\gamma)(t) & = l(t),\quad t>0.
\end{align*}
Both $s_\gamma$ and $r_\gamma$ are nonnegative for all $\gamma\in \iR$. For $\gamma\ge 0$, this is a consequence of the complete positivity of $l$ (see \cite{CN}, \cite{JanI}). If $\gamma<0$
this can be seen, e.g.\ by a simple fixed point argument in the space of nonnegative
$L_1((0,T))$-functions with arbitrary $T>0$ and an appropriate norm. Moreover, $s_\gamma\in H^1_{1,\,loc}([0,\infty))$ for all $\gamma\in \iR$, and if $\gamma\ge 0$, $s_\gamma$ is nonincreasing.

Convolving the $r_\gamma$-equation with $k$ and using that $k\ast l=1$, it follows that $s_\gamma=k\ast r_\gamma$, by uniqueness. Further, we see that
\[
\gamma(1\ast r_\gamma)(t)=1-(k\ast r_\gamma)(t)=1-s_\gamma(t),\quad t>0,
\]
which shows that for $\gamma>0$ the function $r_\gamma$ is integrable on $\iR_+$.

For $\gamma>0$ let $h_\gamma\in L_{1,loc}(\iR_+)$ denote the resolvent kernel associated
with $\gamma l$, that is we have
\begin{equation} \label{hndef}
h_\gamma(t)+\gamma(h_\gamma\ast l)(t)=\gamma l(t),\quad t>0,\,\gamma>0.
\end{equation}
Note that $h_\gamma=\gamma r_\gamma=-\dot{s}_\gamma \in L_{1,\,loc}(\iR_+)$, in particular $h_\gamma$ is nonnegative.
It is well-known that for any $f\in L_p([0,T])$, $1\le p<\infty$, there holds
$h_n\ast f\rightarrow f$ in $L_p([0,T])$ as $n\rightarrow \infty$, see e.g.\ \cite{Za}.

For $\gamma>0$ we set
\begin{equation} \label{kndef}
k_{\gamma}=k\ast h_\gamma.
\end{equation}
It is known (see e.g.\ \cite{Za}) that $k_\gamma=\gamma s_\gamma$, $\gamma>0$, and thus
the kernels $k_\gamma$ are also nonnegative and nonincreasing, and they belong to
$H^1_{1,\,loc}([0,\infty))$ as well.

We next state an important convexity inequality for operators of the form
$\frac{d}{dt}(k\ast \cdot)$. A proof can be found in \cite{KSVZ}.
\begin{lemma} \label{convexFI} Let $T>0$ and $U$ be an open subset of $\iR$. Let further $k\in
H^1_1([0,T])$ be nonnegative and nonincreasing, $H\in C^1(U)$ be convex, $u_0\in U$, and $u\in L_1([0,T])$ with $u(t)\in U$
for a.a. $t\in (0,T)$. Suppose that the functions $H(u)$, $H'(u)u$,
and $H'(u)(\dot{k}\ast u)$ belong to $L_1([0,T])$ (which is the case
if, e.g., $u\in L_\infty([0,T])$). Then
\begin{align} \label{convexfundidentity}
H'(u(t))&\frac{d}{dt}\,\big(k \ast [u-u_0]\big)(t) \ge \;\frac{d}{dt}\,\big(k\ast
[H(u)-H(u_0)]\big)(t),\quad \mbox{a.a.}\;t\in (0,T).
\end{align}
\end{lemma}
\section{Auxiliary results for the linearized problem}
Let $T>0$ and $\Omega\subset \iR^N$ be a bounded domain. In what follows we will use the notation 
$\Omega_T:=(0,T)\times \Omega$ and $\Gamma_T:=(0,T)\times \partial \Omega$. In this section we consider the linear problem
\begin{align}
\partial_t\big(k\ast[v-v_0]\big)-\mbox{div}\big(A(t,x)\nabla v\big) & =m(t,x)v+f(t,x),\quad (t,x)\in \Omega_T,
\nonumber\\
v & =0,\quad (t,x)\in \Gamma_T, \label{linearhilfsproblem}\\
v|_{t=0} & = v_0,\quad x\in \Omega.\nonumber
\end{align}
Here $k$ is a kernel of type $\mathcal{PC}$, $A$ is assumed to satisfy condition ($\mathcal{H}$), $m\in L_\infty(\Omega_T)$, $v_0\in L_2(\Omega)$, and $f\in L_2((0,T);L_2(\Omega))$. Denote by $y_+$ and $y_-:=[-y]_+$ the positive and negative part, respectively, of $y\in \iR$.

We say that $v$ is a {\em weak solution (subsolution, supersolution)} of (\ref{linearhilfsproblem}) if
\begin{itemize}
\item[(a)] $v\in
W(T):=\{w\in  L_2((0,T);H^1_2(\Omega)):\;k\ast w\in C([0,T];L_2(\Omega))\;\mbox{and}\;
(k\ast w)|_{t=0}=0\}$,
\item[(b)]
$v\,\big(v_+,\; v_-\big) \in L_2((0,T);\oH^1_2(\Omega))$, where $\oH^1_2(\Omega):=\overline{C_0^\infty(\Omega)}\,{}^{H^1_2(\Omega)}$,
\item[(c)] for any nonnegative test function
\[
\eta\in H^1_2([0,T];L_2(\Omega))\cap
L_2([0,T];\oH^1_2(\Omega))
\]
with $\eta|_{t=T}=0$ there holds
\[
\int_{0}^{T} \int_\Omega \Big(-\eta_t \big(k\ast [v-v_0]\big)+
(A\nabla v|\nabla \eta)\Big)\,dx\,dt=\,(\le,\,\ge)\,\int_0^T \int_\Omega \big(mv\eta+f\eta\big)\,dx\,dt.
\]
\end{itemize}
Existence and uniqueness of a weak solution to (\ref{linearhilfsproblem})
under the above assumptions
follow from the results in \cite[Theorem 2.1 and Corollary 4.1]{ZWH}.
\begin{satz} \label{linearweakex} 
Let the above assumptions on $\Omega$, $T$, $k$, $A$ and $m$
be fulfilled. Then for any $f\in L_2(\Omega_T)$ and $v_0\in L_2(\Omega)$ the problem (\ref{linearhilfsproblem}) has a unique weak solution $v\in W(T)$ and 
\[
|k\ast v|_{C([0,T];L_2(\Omega))}+|v|_{ L_2((0,T);H^1_2(\Omega))}\le C\big (|v_0|_{L_2(\Omega)}+
|f|_{L_2(\Omega_T)}\big),
\]
where the constant $C$ is independent of $v$, $v_0$, and $f$. 
\end{satz}
Note that $v\in W(T)$ does not entail
$v\in C([0,T];L_2(\Omega))$ in general, so it is not so clear how to interpret the
initial condition. However, once one knows that the functions $v$ and $k\ast
(v-v_0)$ are sufficiently smooth, then $v|_{t=0}=v_0$ is satisfied in an appropriate sense (see \cite{ZWH}). Further, for any weak solution of (\ref{linearhilfsproblem}) we have in addition
$\frac{d}{dt}(k\ast (v-v_0))\in L_2([0,T];H^{-1}_2(\Omega))$, where the time derivative has to be understood in the generalized sense and $H^{-1}_2(\Omega)$
denotes the dual space of $\oH^1_2(\Omega)$, see \cite{ZWH}.

In order to derive a priori estimates for (\ref{linearhilfsproblem}) in a rigorous way one needs a suitable time-regularized
version of (\ref{linearhilfsproblem}). The following equivalent formulation has the advantage that the singular kernel $k$
is replaced by a more regular kernel. For a proof, we refer to \cite[Lemma 3.1]{Za}. The idea behind the regularization
is to replace $\frac{d}{dt}(k\ast \cdot)$ by its Yosida approximations, see also \cite{VZ1}.
\begin{lemma} \label{regWF}
Let the above assumptions be satisfied. Let
$v\in W(T)$ be such that condition (b) above is satisfied. Then $v$ is a weak solution (subsolution, supersolution) of (\ref{linearhilfsproblem})
if and only if for any nonnegative
function $\psi\in \oH^1_2(\Omega)$ there holds
\begin{align} \label{regcond}
\int_\Omega
\Big(\psi\partial_t [k_n\ast(v-v_0)]+\big(h_n\ast[A\nabla v]|\nabla \psi\big)\Big)\,dx=\,(\le,\,\ge)\,
\int_\Omega \psi \big(h_n\ast [mv+f]\big)\,dx,
\end{align}
for a.a.\ $t\in (0,T)$ and all $n\in \iN$. Here $k_n$ is defined as in (\ref{kndef}).
\end{lemma}
We next state the weak maximum principle for (\ref{linearhilfsproblem}) with $m=f=0$. It can be
found in \cite[Theorem 3.1]{Za}.
\begin{satz} \label{MaxPr1} Let the above assumptions on $\Omega$, $T$, $k$, $A$, and $v_0$ be fulfilled.
Assume that $m=f=0$. Then for any weak subsolution (supersolution) $v$ of (\ref{linearhilfsproblem}) there holds for a.a.\
$(t,x)\in \Omega_T$
\[
v(t,x)\le \max\big\{0,\esup_\Omega v_0\big\} \quad\quad \Big(\; v(t,x)\ge \min\big\{0,\einf_\Omega v_0\big\}\;\Big),
\]
provided this maximum (minimum) is finite.
\end{satz}
We also need the following result.
\begin{lemma} \label{fort}
Let the above assumptions on $\Omega$, $T$, $k$, and $A$ be fulfilled. Let $v_0\in L_\infty(\Omega)$,
$f\in L_\infty(\Omega_T)$, and assume $m=0$.
Let $a\in (0,T)$ and $u\in W(a)\cap L_\infty(\Omega_a)$ be the weak solution of  (\ref{linearhilfsproblem}) on $\Omega_a$.
Let $g=f \chi_{(0,a)}(t)$ and suppose that $w\in W(T)$ is the weak solution of the problem
\begin{align*}
\partial_t\big(k\ast[w-v_0]\big)-\mbox{{\em div}}\big(A(t,x)\nabla w\big) & =g(t,x)
,\quad (t,x)\in\Omega_T,
\nonumber\\
w & =0,\quad (t,x)\in \Gamma_T,\\
w|_{t=0} & = v_0,\quad x\in \Omega.\nonumber
\end{align*}
Then there holds for a.a.\
$(t,x)\in \Omega_T$
\[
\min\big\{0,\einf_\Omega v_0, \einf_{\Omega_a} u\big\}\le w(t,x)\le \max\big\{0,\esup_\Omega v_0, \esup_{\Omega_a} u\big\}. 
\]
\end{lemma}
{\em Proof.} We only prove the upper bound for $w$. The lower bound then follows from the upper bound for $-w$. Set 
\[
\kappa:=\max\big\{0,\esup_\Omega v_0, \esup_{\Omega_a} u\big\}.
\]
Testing the time-regularized version of the problem for $t\in (0,T)$ with $w_\kappa^+:=(w-\kappa)_+$ gives
\[
\int_\Omega
\Big(w_\kappa^+\partial_t\big [k_n\ast\big([w-\kappa]-[v_0-\kappa]\big)\big]+\big(h_n\ast[A\nabla w]|\nabla w_\kappa^+\big)\Big)\,dx=
\int_\Omega w_\kappa^+ \big(h_n\ast g\big)\,dx,
\]
for a.a.\ $t\in (0,T)$. By Lemma \ref{convexFI} with $H(y)=\frac{1}{2}(y_+)^2$, $y\in \iR$, applied to the first term, and
since $H(v_0-\kappa)=0$, it follows that
\[
\int_\Omega
\Big(\partial_t\big [k_n\ast H(w-\kappa)\big]+\big(h_n\ast[A\nabla w]|\nabla w_\kappa^+\big)\Big)\,dx\le
\int_\Omega w_\kappa^+ \big(h_n\ast g\big)\,dx,
\]
for a.a.\ $t\in (0,T)$. Convolving this inequality with the nonnegative kernel $l$, using that
\[
l\ast \partial_t\big(k_n\ast H(w-\kappa)\big)=\partial_t \big(l\ast k\ast h_n\ast H(w-\kappa)\big)
=h_n\ast H(w-\kappa)
\]
(note that $(k_n\ast H(w-\kappa))|_{t=0}=0$), and sending $n\to \infty$, we infer that
\[
|w_\kappa^+(t,\cdot)|_{L_2(\Omega)}^2+2 l\ast \int_\Omega \big(A\nabla w|\nabla [w_\kappa^+]\big)\,dx\le 2 l\ast\int_\Omega g w_\kappa^+\,dx,\quad\mbox{a.a.}\;t\in (0,T).
\]
But now $g w_\kappa^+=0$ a.e.\ in $\Omega_T$, since on $\Omega_a$ $w$ and $u$ coincide (by uniqueness) and thus
$w_\kappa^+=0$ a.e.\ in $\Omega_a$, by definition of $\kappa$. Using this and ($\mathcal{H}$), it follows that
\[
|w_\kappa^+(t,\cdot)|_{L_2(\Omega)}^2\le 0,\quad\mbox{a.a.}\;t\in (0,T),
\]
that is
$w\le \kappa$ a.e.\ in $\Omega_T$. \hfill $\square$

${}$

\noindent The next result provides the comparison principle for (\ref{linearhilfsproblem}).
\begin{satz} \label{comp}
Let the above assumptions on $\Omega$, $T$, $k$, $A$ and $m$
be fulfilled. Let $f\in L_2(\Omega_T)$ and $v_0\in L_2(\Omega)$. Suppose that $u\in W(T)$ is a weak subsolution of (\ref{linearhilfsproblem}) and that $v\in W(T)$ is a weak supersolution of (\ref{linearhilfsproblem}). Then $u\le v$ a.e.\ in $\Omega_T$.
\end{satz}
{\em Proof.} We give only a sketch of the proof. Using Lemma \ref{regWF} and setting $w=u-v$ we have $w\in W(T)$, $w_+
\in L_2((0,T);\oH^1_2(\Omega))$, and for any nonnegative
function $\psi\in \oH^1_2(\Omega)$ there holds
\begin{align} \label{comp1}
\int_\Omega
\Big(\psi\partial_t [k_n\ast w]+\big(h_n\ast[A\nabla w]|\nabla \psi\big)\Big)\,dx\le
\int_\Omega \psi \big(h_n\ast [mw]\big)\,dx,
\end{align}
for a.a.\ $t\in (0,T)$ and all $n\in \iN$. We test this inequality with $w_+$ and proceed similarly
as in the proof of Lemma \ref{fort}, that is, we apply the convexity inequality from Lemma \ref{convexFI},
convolve the resulting inequality with $l$ and let $n\to \infty$. Using also assumption ($\mathcal{H}$)
this leads to the inequality 
\[
|w_+(t)|^2_{L_2(\Omega)}\le 2|m|_{L_\infty(\Omega_T)}\big( l\ast |w_+(\cdot)|^2_{L_2(\Omega)}\big)(t), \quad\mbox{a.a.}\;t\in (0,T).
\]
Since $l$ is nonnegative, this inequality implies that $|w_+(t)|_{L_2(\Omega)}^2=0$ a.e.\ in $(0,T)$, by the abstract Gronwall lemma \cite[Prop.
7.15]{ZeidlerI}, i.e. $u\le v$ a.e.\ in $\Omega_T$. \hfill $\square$

${}$

\noindent By means of the comparison principle we obtain the following result.
\begin{korollar} \label{compfbounded}
Let the above assumptions on $\Omega$, $T$, $k$, and $A$ be satisfied. Let $f\in L_\infty(\Omega_T)$ and $v_0\in L_\infty(\Omega)$ and assume that $m=0$. Let $v\in W(T)$ be the weak solution of (\ref{linearhilfsproblem}). Then
$v\in L_\infty(\Omega_T)$ and
\begin{equation} \label{linearbound}
|v|_{L_\infty(\Omega_T)}\le |v_0|_{L_\infty(\Omega)}+(1\ast l)(T) |f|_{L_\infty(\Omega_T)}.
\end{equation}
\end{korollar}
{\em Proof.} Setting $M=|f|_{L_\infty(\Omega_T)}$  we have 
\begin{align*}
\partial_t\big(k\ast[v-v_0]\big)-\mbox{div}\big(A(t,x)\nabla v\big) & \le M,\quad (t,x)\in\Omega_T,
\nonumber\\
v & =0,\quad (t,x)\in \Gamma_T, \\
v|_{t=0} & =v_0,\quad x\in \Omega,\nonumber
\end{align*}
in the weak subsolution sense. On the other hand, the function $w(t,x)=w_0+M(1\ast l)(t)$
with $w_0=|v_0|_{L_\infty(\Omega)}$ satisfies
 \begin{align*}
\partial_t\big(k\ast[w-w_0]\big)-\mbox{div}\big(A(t,x)\nabla w\big) & = M,\quad (t,x)\in \Omega_T,\nonumber\\
w & \ge 0,\quad (t,x)\in \Gamma_T, \\
w|_{t=0} & = w_0,\quad x\in \Omega\nonumber.
\end{align*}
Since $v_0\le w_0$, the comparison principle implies $v\le w$ a.e.\ in
$\Omega_T$. Replacing $v$ with $-v$ and $v_0$ with $-v_0$ the same
argument shows that $-v\le w$. Hence $|v|\le w$ a.e.\ in
$\Omega_T$. \hfill $\square$

${}$

\noindent The subsequent lemma says that the positive part of a bounded weak subsolution to 
(\ref{linearhilfsproblem}) with $m=0$ is a bounded weak subsolution to a related problem. By $\chi_M$ we mean the 
characteristic function of the set $M$.
\begin{lemma} \label{pospartsub}
Let the above assumptions on $\Omega$, $T$, $k$, and $A$ be satisfied. Let $f\in L_\infty(\Omega_T)$ and $v_0\in L_\infty(\Omega)$ and assume that $m=0$. Let $v\in W(T)$ be a bounded weak subsolution of (\ref{linearhilfsproblem}). Then the positive part of $v$ is a bounded weak subsolution of the problem
\begin{align}
\partial_t\big(k\ast[w-w_0]\big)-\mbox{{\em div}}\big(A(t,x)\nabla w\big) & =f \chi_{\{v\ge 0\}},\quad (t,x)\in \Omega_T,
\nonumber\\
w & =0,\quad (t,x)\in \Gamma_T, \label{pospartprob}\\
w|_{t=0} & = w_0,\quad x\in \Omega,\nonumber
\end{align}
where $w_0=(v_0)_+$. 
\end{lemma}
{\em Proof.} Note first that $v\in W(T)\cap L_\infty(\Omega_T)$ implies that $v_+$ belongs to the same space. The claimed subsolution property of $v_+$ can be shown by the same line of arguments as in \cite[Section 4]{VZ}. The idea is to test the time-regularized subsolution inequality for $v$ with a suitable regularization $H_\varepsilon'(v)$ of 
$\chi_{\{v\ge 0\}}$ and to apply Lemma \ref{convexFI} to the convex function $H_\varepsilon$. Letting finally $\varepsilon\to 0$ yields the assertion. \hfill $\square$

${}$

\noindent The last result of this section provides, among others, sufficient conditions for the stability of the zero function for the linear problem (\ref{linearhilfsproblem}) with $f=0$. This result will also be crucial
for the nonlinear stability analysis. Here we will assume in addition that
$A$ is independent of $t$. By $\lambda_1>0$ we mean the first eigenvalue
of the negative Dirichlet-Laplacian $(-\Delta_D)$ in $L_2(\Omega)$.
If $A$ is also symmetric, by $\lambda_*>0$ we denote the smallest eigenvalue of the operator
$\mathcal{L}v=-\mbox{div}\big(A(x)\nabla v\big)$ (with Dirichlet boundary condition) in $L_2(\Omega)$. 
\begin{satz} \label{linearestimate}
Let $T>0$ and $\Omega$ be a bounded domain in $\iR^N$. Suppose that
($\mathcal{H}$) is satisfied and that $A$ is independent of $t$.
Suppose that $(k,l)\in \mathcal{PC}$ for some
$l\in L_{1,\,loc}(\iR_+)$. Let $v_0\in L_\infty(\Omega)$ and $m,\,f\in L_\infty(\Omega_T)$. 
Let $v\in W(T)$ be the weak solution of (\ref{linearhilfsproblem}). Then
$v\in L_\infty(\Omega_T)$ and there exists a constant ${C}>0$ that is
independent of $v, v_0, f$ such that
\begin{equation} \label{linbound2}
|v|_{L_\infty(\Omega_T)}\le {C}\big(|v_0|_{L_\infty(\Omega)}+|f|_{L_\infty(\Omega_T)}\big).
\end{equation}
Assume in addition that one of the two following {\em stability conditions}
is satisfied.
\begin{itemize}
\item [(a)] $A$ is also symmetric and $\esup_{(t,x)\in \Omega_T} m(t,x)<\lambda_*$.
\item [(b)] $\esup_{(t,x)\in \Omega_T} m(t,x)<\nu \lambda_1$.
\end{itemize}
Then the constant ${C}$ in (\ref{linbound2}) can be chosen independent of $T$ and in the special case $f=0$ there holds
\begin{equation} \label{linboundextra}
|v(t,x)|\le \tilde{C} s_{\delta_0}(t)|v_0|_{L_\infty(\Omega)},\quad \mbox{a.a.}\,(t,x)\in\Omega_T,
\end{equation}
with some $\delta_0>0$ and $\tilde{C}>0$ independent of $v$, $v_0$ and $T$. Here $s_{\delta_0}$ is the relaxation function whose definition was given in Section \ref{Prelim}. In particular, if condition a) resp.\ b) holds for all $T>0$ then
the solution $v$ of (\ref{linearhilfsproblem}) with $T=\infty$ tends to $0$ as $t\to \infty$ whenever $l\notin L_1(\iR_+)$.
\end{satz}
{\em Proof.} We give the argument in case that $A$ is also symmetric. The proof in the nonsymmetric
case is the same, one only has to replace $\lambda_*$ by $\nu \lambda_1$ in all formulas.

To establish the boundedness of $v$, we introduce the function $\vartheta\in H^1_2(\Omega)$ as the
weak solution of the elliptic problem
\begin{align*}
-\mbox{div}\big(A(x)\nabla \vartheta\big) & = (\lambda_*-\varepsilon)\vartheta,\quad x\in \Omega,\\
\vartheta & = 1,\quad x\in \partial\Omega,
\end{align*}
where we fixed some $\varepsilon\in (0,\lambda_*)$.
Note that $\vartheta$ is well-defined, thanks to the assumptions on $A$ and the Lax-Milgram lemma.
The boundary condition on $\vartheta$ has to be interpreted in the weak sense as $\vartheta-1\in \oH^1_2(\Omega)$.
The comparison principle implies that $\vartheta\ge 1$ in $\Omega$, and by elliptic regularity
theory we also have $\vartheta\in L_\infty(\Omega)$.

Next, set $\mu:=\mbox{ess}\,\sup_{(t,x)\in \Omega_T}m(t,x)$ and $\kappa:=\mu-\lambda_*+\varepsilon$ and consider the positive and bounded function
\[
\zeta(t,x)=s_{-\kappa}(t)\vartheta(x)|v_0|_\infty+(1\ast r_{-\kappa})(t)\vartheta(x)|f|_\infty,\quad (t,x)\in \Omega_T.
\]
By the definition of $s_\gamma$ and $r_\gamma$ (cf.\ Section \ref{Prelim}) we have
\begin{align*}
\partial_t\big(k\ast[s_\gamma-1]\big) & =-\gamma s_\gamma,\quad t>0,\\
\partial_t\big(k\ast[1\ast r_\gamma]\big) & =-\gamma \,1\ast r_\gamma+1,\quad t>0,
\end{align*}
and thus with $\zeta_0(x):=\vartheta(x)|v_0|_\infty$,
\begin{align*}
\partial_t\big(k\ast[\zeta-\zeta_0]\big)-\mbox{div}\big(A(x)\nabla \zeta\big) = &\; \kappa \zeta+\vartheta(x)|f|_\infty+ (\lambda_*-\varepsilon)\zeta\\
= &\;\mu \zeta+\vartheta(x)|f|_\infty\\
\ge &\; m(t,x)\zeta+f(t,x).
\end{align*}
Since also $\zeta_0\ge |v_0|_\infty\ge v_0$ in $\Omega$, $\zeta$ is a supersolution of (\ref{linearhilfsproblem}), which implies $v\le \zeta$, by the comparison principle, Theorem \ref{comp}.

Looking at $-v$ instead of $v$, the above argument shows that $v\ge -\zeta$ in $\Omega_T$.
Hence $|v|\le \zeta$, which in turn yields boundedness of $v$, together with (\ref{linbound2}).

Finally, suppose that the stability condition $\mu<\lambda_*$ holds. Then we may select $\varepsilon
\in (0,\lambda_*)$ such that $\kappa=\mu-\lambda_*+\varepsilon<0$. In this case $s_{-\kappa}$ is
nonincreasing and $r_{-\kappa}\in L_1(\iR_+)$, cf.\ Section \ref{Prelim}. Therefore
\[
|v|_{L_\infty((0,T)\times \Omega)} \le |\vartheta|_\infty\big(|v_0|_\infty+ |r_{-\kappa}|_{L_1(\iR_+)}|f|_\infty\big),
\]
which proves the statement on the constant $T$. Assertion (\ref{linboundextra}) follows from $|v|\le \zeta$ on $\Omega_T$ and the structure of $\zeta$ with $\delta_0=-\kappa>0$. As to the last claim we refer to \cite[Lemma 6.1]{VZ}.
 \hfill $\square$

${}$

\noindent To illustrate our linear stability result, we give some examples of pairs $(k,l)\in {\cal PC}$
and discuss the decay behaviour of the corresponding relaxation function $s_\mu$ for $\mu>0$.
These and further examples can be found in \cite[Section 6]{VZ}.
\begin{bei} \label{kerne}
{\em 
a) The classical time fractional case.
We consider the pair
\begin{equation} \label{Fbeifrac}
(k,l)=(g_{1-\alpha},g_\alpha),\quad \mbox{where}\;\alpha\in (0,1).
\end{equation}
In this case
\[
s_\mu(t)=E_\alpha(-\mu t^\alpha),\quad\mbox{where}\;E_\alpha(z):=\sum_{j=0}^\infty
\,\frac{z^j}{\Gamma(\alpha j+1)}\,,\;z\in \iC,
\]
is the well-known Mittag-Leffler function (see e.g.\ \cite{KST}), which satisfies the estimate
\[
\frac{1}{1 + \Gamma (1-\alpha) x} \leq E_{\alpha}(-x) \leq \frac{1}{1+\frac{x}{\Gamma(1+\alpha)}},\,\quad x\ge 0,
\]
see \cite[Example 6.1]{VZ}. Thus with $C(\alpha)=\Gamma(1+\alpha)^{-1}$ we obtain for $\mu>0$ the algebraic decay estimate
\[
s_\mu(t)\le \frac{1}{1+C(\alpha) \mu t^\alpha},\quad t\ge 0.
\]

b) The time fractional case with exponential weight. We consider
\[
k(t)=g_{1-\alpha}(t)e^{-\gamma t},\quad l(t)=g_{\alpha}(t)e^{-\gamma
t}+\gamma(1\ast[g_{\alpha}e^{-\gamma\cdot}])(t), \quad t>0,
\]
with $\alpha\in (0,1)$ and $\gamma>0$. Let $\mu>0$ be fixed. Then 
$s_\mu(t)\le Me^{-\omega t}$ for all $t\ge 0$
where $M$ is independent of $t$ and $\omega\in (0,\gamma)$ is the unique solution of $\omega=\mu(\gamma-\omega)^{1-\alpha}$, see \cite[Example 6.2]{VZ}.

c) An example of ultraslow diffusion. We consider the pair (\ref{ultrapair}) already mentioned in the introduction, that is
\[
k(t)=\int_0^1 g_\beta(t)\,d\beta,\quad l(t)=\int_0^\infty \,\frac{e^{-st}}{1+s}\,{ds},\quad t>0.
\]
It is shown in \cite[Example 6.5]{VZ} that $(k,l)\in \mathcal{PC}$ and that there is a number $T>1$ independent of $\mu\ge 0$ such that
\[
s_\mu(t)\le \,\frac{1}{1+\,\frac{\mu}{2}\,\log t},\quad t\ge T.
\]

}
\end{bei}
\section{Well-posedness}
We have the following result on the well-posedness of the nonlinear problem (\ref{subproblem}).
\begin{satz} \label{wellposedsubdiffusion}
Let $\Omega$ be a bounded domain in $\iR^N$, and $f\in C^{1-}(I)$, where $I\neq \emptyset$ is an open interval in $\iR$. Suppose that
($\mathcal{H}$) is satisfied and that $(k,l)\in \mathcal{PC}$ for some
$l\in L_{1,\,loc}(\iR_+)$. Suppose further that $u_0\in L_\infty(\Omega)$ with $I_0:=[\einf_\Omega u_0,\esup_\Omega u_0]\subset I$ and that $0\in I$. Then the following holds.

(i) There exists a maximal existence time $t_*\in (0,\infty]$ such that problem (\ref{subproblem}) admits for any $a\in (0,t_*)$ a unique solution
\[
u\in Z(a):=L_\infty((0,a)\times \Omega)\cap L_2((0,a);\oH^1_2(\Omega)).
\]

(ii) Given $\varepsilon>0$ and $\delta\in (0,\varepsilon)$, there exists $a\in (0,t_*)$ such that
we have the implication
\[
|u_0|_{L_\infty(\Omega)}\le \delta \;\Rightarrow\;|u|_{L_\infty((0,a)\times \Omega)}
\le \varepsilon.
\]

(iii) If in addition $u_0\ge 0$ a.e.\ in $\Omega$ and $f(0)\ge 0$ then the solution
$u$ is nonnegative a.e.\ in $(0,t_*)\times \Omega$.

(iv) If $I=\iR$ and $t_*<\infty$ then
\[
t_*=\sup\{a>0:\,|u|_{L_\infty((0,a)\times \Omega)}<\infty\}.
\]

(v) Let $a\in (0,t_*)$ and suppose that $\bar{u}\in Z(a)$ solves  (\ref{subproblem}) on
$\Omega_a$ with $u_0$ replaced with 
$\bar{u}_0\in L_\infty(\Omega)$ satisfying $[\einf_\Omega \bar{u}_0,\esup_\Omega \bar{u}_0]\subset I$.
Let $K$ be the union of the (essential) ranges
$K=u\big(\Omega_a\big)\cup \bar{u}\big(\Omega_a\big)\subset I.$ 
Then there holds the stability estimate
\begin{equation} \label{stabinitial}
|u-\bar{u}|_{L_\infty(\Omega_a)}\le C|u_0-\bar{u}_0|_{L_\infty(\Omega)},
\end{equation}
where the constant $C=C(f,K,a,l)$.
\end{satz}

${}$

\noindent In the following, by a {\em global solution} of (\ref{subproblem}) we mean a solution with $t_*=\infty$.

${}$

\noindent {\em Proof of Theorem \ref{wellposedsubdiffusion}.} {\bf 1. Stability estimate w.r.t.\ initial value and uniqueness.} Let $a>0$ and set $X(a):=L_\infty(\Omega_a)$. Suppose that $u\in Z(a)$ solves  (\ref{subproblem}) on
$\Omega_a$ and that $\bar{u}\in Z(a)$ solves  (\ref{subproblem}) with $u_0$ replaced with $\bar{u}_0$ on
$\Omega_a$. Setting $v=u-\bar{u}$ and $v_0=u_0-\bar{u}_0$ we then have (in the weak sense)
\begin{align*}
\partial_t\big(k\ast [v-v_0]\big)-\mbox{div}\big(A(t,x)\nabla v\big) & =f(u)-f(\bar{u})
,\quad (t,x)\in \Omega_a,
\nonumber\\
v & =0,\quad (t,x)\in \Gamma_a,\\
v|_{t=0} & = v_0,\quad x\in \Omega.\nonumber
\end{align*}
Since $u$ and $\bar{u}$ are (essentially bounded) solutions, we can modify $u$ and $\bar{u}$ on a set 
${\mathcal N}\subset \Omega_a$ of measure zero such that the union of their ranges
$
K:=u\big(\Omega_a\big)\cup \bar{u}\big(\Omega_a\big)
$ 
is a compact subset of $I$. Since $f\in C^{1-}(I)$, $f$ is Lipschitz continuous on $K$. Denoting by $L$ the corresponding Lipschitz constant we have
\[
|f(u)-f(\bar{u})|\le L|u-\bar{u}|=L|v|\quad \mbox{in}\;\Omega_a,
 \] 
and thus (in the sense of a weak subsolution)
\begin{align*}
\partial_t\big(k\ast [v-v_0]\big)-\mbox{div}\big(A(t,x)\nabla v\big) & \le L|v|
,\quad (t,x)\in \Omega_a,
\nonumber\\
v & =0,\quad (t,x)\in \Gamma_a,\\
v|_{t=0} & = v_0,\quad x\in \Omega.\nonumber
\end{align*}

Since $|v| \chi_{\{v\ge 0\}}=v_+$, it follows from Lemma \ref{pospartsub} that
\begin{align*}
\partial_t\big(k\ast [v_+-(v_0)_+]\big)-\mbox{div}\big(A(t,x)\nabla v_+\big) & \le Lv_+
,\quad (t,x)\in \Omega_a,
\nonumber\\
v_+ & =0,\quad (t,x)\in \Gamma_a,\\
v_+|_{t=0} & = (v_0)_+,\quad x\in \Omega,\nonumber
\end{align*}
in the weak sense.

Setting $\zeta_0=|(v_0)_+|_{L_\infty(\Omega)}$, the nonnegative function
$\zeta(t)=s_{-L}(t)\zeta_0$ (see Section \ref{Prelim} for the definition of the relaxation function $s_{-L}$)
satisfies
\begin{align*}
\partial_t\big(k\ast [\zeta-(v_0)_+]\big)-\mbox{div}\big(A(t,x)\nabla \zeta\big) & \ge L\zeta
,\quad (t,x)\in \Omega_a,
\nonumber\\
\zeta & \ge 0,\quad (t,x)\in \Gamma_a,\\
\zeta|_{t=0} & \ge  (v_0)_+,\quad x\in \Omega,\nonumber
\end{align*}
in the weak sense. Thus the comparison principle, Theorem \ref{comp}, implies that
\[
\big(u-\bar{u}\big)_+=v_+\le \zeta(t)=s_{-L}(t)|(u_0-\bar{u}_0)_+|_{L_\infty(\Omega)}\quad
\mbox{in}\;\Omega_a.
\]
Analogously, one obtains a corresponding upper bound for $\big(\bar{u}-{u}\big)_+$. Combining
both estimates gives
\begin{equation} \label{vzeta}
|u-\bar{u}|_{L_\infty(\Omega_a)}\le C|u_0-\bar{u}_0|_{L_\infty(\Omega)},
\end{equation}
where $C=C(L,a,l)$. This shows the stability estimate (\ref{stabinitial}). In particular, taking
$u_0=\bar{u}_0$ we obtain uniqueness for problem (\ref{subproblem}). 

{\bf 2. Local existence.} Let $T>0$ and $w\in Z(T)$ be the weak solution of the linear problem
\begin{align*}
\partial_t\big(k\ast[w-u_0]\big)-\mbox{div}\big(A(t,x)\nabla w\big) & =0
,\quad (t,x)\in\Omega_T,
\nonumber\\
w & =0,\quad (t,x)\in \Gamma_T,\\
w|_{t=0} & = u_0,\quad x\in \Omega.\nonumber
\end{align*}
The function $w$ is well-defined, thanks to Theorem \ref{linearweakex} and Theorem \ref{MaxPr1}. Moreover, 
\begin{equation} \label{wbounds}
\min\big\{0,\einf_\Omega u_0\big\} \le w(t,x)
\le \max\big\{0,\esup_\Omega u_0\big\},
\end{equation}
for a.a.\ $(t,x)\in \Omega_T$, by the maximum principle.
Recall the assumptions $I_0\subset I$ and $0\in I$. 
So, in view of (\ref{wbounds}) the essential range of $w$ is contained in a compact subset of the open interval $I$. 

We next fix $\rho>0$ such that $I_1:=[\einf_{\Omega_T} w-\rho,\esup_{\Omega_T} w+\rho]\subset I$.
For $a\in (0,T]$  
we introduce the non-empty set
\[
\Sigma(a,\rho):=\{v\in X(a):\,|v-w|_{X(a)}\le \rho\}.
\]
Invoking Corollary \ref{compfbounded}, we define a map $\Phi:\Sigma(a,\rho)\rightarrow Z(a)$ by assigning to $u$ the weak solution $v=\Phi(u)$ of the linear problem
\begin{align*}
\partial_t\big(k\ast[v-u_0]\big)-\mbox{div}\big(A(t,x)\nabla v\big) & =
f(u),\quad (t,x)\in\Omega_a,
\nonumber\\
v & =0,\quad (t,x)\in \Gamma_a,\\
v|_{t=0} & = u_0,\quad x\in \Omega.\nonumber
\end{align*}
Note that in view of the choice of $\rho$, the term $f(u)$ is well-defined for any $u\in\Sigma(a,\rho)$. 
We will show that for sufficiently small $a$ the map
$\Phi$ leaves $\Sigma(a,\rho)$ invariant and becomes a strict contraction in
$X(a)$.

Let $u\in\Sigma(a,\rho)$ and $v=\Phi(u)$. Then the difference $v-w$ solves
the problem
\begin{align*}
\partial_t\big(k\ast[v-w]\big)-\mbox{div}\big(A(t,x)\nabla (v-w)\big) & =
f(u)\;\; \mbox{in}\;\Omega_a,
\nonumber\\
v-w & =0\;\;\mbox{on}\;\Gamma_a,\\
(v-w)|_{t=0} & = 0\;\;\mbox{in}\;\Omega.\nonumber
\end{align*}
Corollary \ref{compfbounded} yields the estimate
\begin{align} \label{invar}
|v-w|_{X(a)} & \le \delta(a)|f(u)|_{X(a)}\nonumber\\
& \le \delta(a)\big(|f(u)-f(w)|_{X(a)}+|f(w)|_{X(a)}\big),
\end{align}
where the constant $\delta(a)=(1\ast l)(a)\to 0$ as $a\to 0$.

Let $L$ be the Lipschitz constant of $f$ on the interval $I_1$.
Since $u$ and $w$ take values in $I_1$,  (\ref{invar}) and the Lipschitz estimate for $f$ imply
\begin{align}
|v-w|_{X(a)} &\, \le \delta(a)\big( L|u-w|_{X(a)}+|f(w)|_{X(T)}\big)\nonumber\\
& \, \le  \delta(a)\big(L\rho+|f(w)|_{X(T)}\big). \label{selfmap}
\end{align}

Next, let $u_1,\,u_2\in\Sigma(a,\rho)$ and $v_i=\Phi(u_i)$, $i=1,2$. Then $v_1-v_2$ solves
the problem
\begin{align*}
\partial_t\big(k\ast[v_1-v_2]\big)-\mbox{div}\big(A(t,x)\nabla (v_1-v_2)\big) & =
f(u_1)-f(u_2)\quad \mbox{in}\;\Omega_a,
\nonumber\\
v_1-v_2 & =0\;\;\mbox{on}\;\Gamma_a,\\
(v_1-v_2)|_{t=0} & = 0\;\;\mbox{in}\;\Omega,\nonumber
\end{align*}
and thus
\begin{equation} \label{contract}
|v_1-v_2|_{X(a)}\le \delta(a)|f(u_1)-f(u_2)|_{X(a)}\le \delta(a) L |u_1-u_2|_{X(a)}.
\end{equation}

Choosing $a$ so small that
\[
\delta(a)(L\rho+|f(w)|_{X(T)})\le \rho\quad \mbox{and}\quad
\delta(a)L\le \frac{1}{2},
\]
we see from (\ref{selfmap}) and (\ref{contract}) that we may apply the contraction mapping principle
to $\Phi$. The unique fixed point of $\Phi$ in the set $\Sigma(a,\rho)$ lies in $Z(a)$ and is a local in time weak
solution of (\ref{subproblem}).

{\bf 3. The maximally defined solution.} The local solution $u\in Z(a)$ obtained in the
second part can be extended to some larger time interval $(0,a+a_1)$. In fact, let $T>a$ and
define now the reference function $w$ as solution of the linear problem
\begin{align*}
\partial_t\big(k\ast[w-u_0]\big)-\mbox{div}\big(A(t,x)\nabla w\big) & =g(t,x)
\quad \mbox{in}\;\Omega_T,
\nonumber\\
w & =0\quad \mbox{on}\;\Gamma_T,\\
w|_{t=0} & = u_0\quad \mbox{in}\; \Omega,\nonumber
\end{align*}
where $g(t,x)=f(u(t,x)) \chi_{(0,a)}(t)$.
Note that $w|_{\Omega_a}=u$, by uniqueness. By Lemma \ref{fort}, 
\[
\min\big\{0,\einf_\Omega u_0, \einf_{\Omega_a} u\big\}\le w(t,x)\le \max\big\{0,\esup_\Omega u_0, \esup_{\Omega_a} u\big\},
\]
for a.a.\ $(t,x)\in \Omega_T$.

Next, fix $\rho>0$ such that $I_2:=[\einf_{\Omega_T} w-\rho,\esup_{\Omega_T} w+\rho]
\subset I$. For $a_1\in (0,T-a]$ we introduce the set
\[
\Sigma(a,a_1,\rho):=\{v\in X(a+a_1):\,v|_{\Omega_a}=u\;\mbox{a.e. in}\,\Omega_a,\;|v-w|_{X(a+a_1)}\le \rho\},
\]
which contains $w$. Define the mapping $\Phi:\Sigma(a,a_1,\rho)\rightarrow
Z(a+a_1)$, which assigns to $\bar{u}\in \Sigma(a,a_1,\rho)$ the solution
$v=\Phi(\bar{u})$ of the linear problem
\begin{align*}
\partial_t\big(k\ast[v-u_0]\big)-\mbox{div}\big(A(t,x)\nabla v\big) & =
f(\bar{u})\quad \mbox{in}\;\Omega_{a+a_1},
\nonumber\\
v & =0\quad \mbox{on}\;\Gamma_{a+a_1},\\
v|_{t=0} & = u_0\quad \mbox{in}\;\Omega.\nonumber
\end{align*}
Since $\bar{u}|_{\Omega_a}=u$, we have $v|_{\Omega_a}=u$, by uniqueness.

Setting $z:=v-w$, it is evident that $z|_{\Omega_a}=0$. For $t\in (a,a+a_1)$ we shift the time by setting
$s=t-a$ and $\tilde{z}(s,x)=z(t,x)$ as well as $\tilde{A}(s,x)=A(t,x)$. Since $z|_{\Omega_a}=0$ we have
\begin{align*}
(k\ast z)(t,x)=\int_a^t k(t-\tau)z(\tau,x)\,d\tau=\int_0^{t-a} k(t-a-\sigma) z(\sigma+a,x)\,d\sigma
\end{align*}
and thus
\[
\partial_t(k\ast z)(t,x)=\partial_s(k\ast \tilde{z})(s,x).
\]
Consequently, the problem for $\tilde{z}$ (to be understood in the weak sense) then reads as
\begin{align*}
\partial_s\big(k\ast \tilde{z}\big)-\mbox{div}\big(\tilde{A}(s,x)\nabla \tilde{z}\big) & =
f(\bar{u}(s+a,x)),\quad (s,x)\in \Omega_{a_1},
\nonumber\\
\tilde{z} & =0,\quad (s,x) \in \Gamma_{a_1},\\
\tilde{z}|_{s=0} & = 0,\quad x\in \Omega.\nonumber
\end{align*}
By Corollary \ref{compfbounded}, it follows that
\[
|\tilde{z}|_{X(a_1)}\le \delta(a_1)|f(\bar{u}(a+\cdot,\cdot))|_{X(a_1)}
\]
with $\delta(a_1)=(1\ast l)(a_1)\to 0$ as $a_1\to 0$. Denoting by $L$ the Lipschitz constant of $f$
on the interval $I_2$ we may argue as in (\ref{invar}), (\ref{selfmap}) to get
\[
|v-w|_{X(a+a_1)}\le \delta(a_1)\big(L\rho+|f(w)|_{X(T)}\big).
\]

Using the same time-shifting trick, we may repeat the argument from Step 2 for the
contraction estimate to see that for any $u_1,u_2\in \Sigma(a,a_1,\rho)$ and
$v_i=\Phi(u_i)$, $i=1,2$,
\[
|v_1-v_2|_{X(a+a_1)}\le \delta(a_1)L
|u_1-u_2|_{X(a+a_1)}.
\]
We see that for sufficiently small $a_1$ the contraction principle applies, yielding
a unique fixed point of $\Phi$ in $\Sigma(a,a_1,\rho)$, which is the unique weak
solution of (\ref{subproblem}) on $\Omega_{a+a_1}$.

Repeating this argument we obtain a maximal interval of existence $(0,t_*)$ with
$t_*\in (0,\infty]$ (recall that $T>0$ was arbitrarily fixed) that is the
supremum of all $\tau>0$ such that (\ref{subproblem}) has a unique solution
$u\in Z(\tau)$. This proves (i).

{\bf 4. Proof of (ii).} Fix $a_0\in (0,t_*)$ and set $m_1=\einf_{\Omega_{a_0}} u$ and
$m_2=\esup_{\Omega_{a_0}} u$. Putting $g=f(u)$ on $\Omega_{a_0}$ we have 
$|g|_{X(a_0)}\le |f|_{L_\infty([m_1,m_2])}=:M$ and for any $a\in (0,a_0]$ we have
 \begin{align*}
\partial_t\big(k\ast[u-u_0]\big)-\mbox{div}\big(A(t,x)\nabla u\big) & = g,\quad (t,x)\in \Omega_a,
\nonumber\\
u & =0,\quad (t,x)\in \Gamma_a, \\
u|_{t=0} & = u_0,\quad x\in \Omega,\nonumber
\end{align*}
in the weak sense. By Corollary \ref{compfbounded},
\[
|u|_{X(a)}\le |u_0|_{L_\infty(\Omega)}+M(1\ast l)(a) \to |u_0|_{L_\infty(\Omega)}
\quad \mbox{as}\; a\to 0.
\]

{\bf 5. Proof of (iii).} Suppose $u_0$ is nonnegative and that $f(0)\ge 0$. Let $u$ be the maximally defined solution of (\ref{subproblem}) on $\Omega_{t_*}$.
Let $a\in (0,t_*)$ and set $m_1=\einf_{\Omega_{a}} u$ and
$m_2=\esup_{\Omega_{a}} u$. Let $L$ be the Lipschitz constant of
$f$ on the interval $[\min\{0,m_1\},\max\{0,m_2\}]\subset I$. Then
\[
f(u)=-\big(f(0)-f(u)\big)+f(0)\ge -L|u|\quad \mbox{a.e. in}\;\Omega_a.
\]
This implies that $v:=-u$ satisfies (with $v_0=-u_0$)
\begin{align*}
\partial_t\big(k\ast[v-v_0]\big)-\mbox{div}\big(A(t,x)\nabla v\big) & \le L|v|,\quad (t,x)\in \Omega_a,
\nonumber\\
v & =0,\quad (t,x)\in \Gamma_a, \\
v|_{t=0} & = v_0,\quad x\in \Omega,\nonumber
\end{align*}
in the weak sense. We can now argue as in Step 1 to obtain that
\begin{align*}
\partial_t\big(k\ast [v_+-(v_0)_+]\big)-\mbox{div}\big(A(t,x)\nabla v_+\big) & \le Lv_+
,\quad (t,x)\in \Omega_a,
\nonumber\\
v_+ & =0,\quad (t,x)\in \Gamma_a,\\
v_+|_{t=0} & = (v_0)_+=0,\quad x\in \Omega,\nonumber
\end{align*}
which in turn implies $v_+=0$ in $\Omega_a$, by the same comparison argument as in Step 1 and
since $(v_0)_+=0$.  This shows nonnegativity of $u$ a.e.\ in
$\Omega_a$. Since $a\in (0,t_*)$ was arbitrary, this proves claim (iii).

{\bf 6. Proof of (iv).} Let $I=\iR$ and assume that $t_*<\infty$. Suppose that there is $b>0$ such that
$|u|_{X(a)}\le b$ for all $a\in (0,t_*)$. We want to show that this contradicts the definition of $t_*$. 

We follow the line of arguments given in Step 3. We may take $a<t_*$ with $t_*-a$ as small as we want.
By the uniform bound for $|u|_{X(a)}$ the Lipschitz constant $L$ in Step 3 can be chosen independently
of $a\in (0,t_*)$, and thus also $a_1$ can be selected independently of the size of $t_*-a$. This means
that for $t_*-a$ sufficiently small the number $a+a_1$ exceeds $t_*$, that is, the solution can be extended
to some interval $[0,t_*+\epsilon]$ with $\epsilon>0$, a contradiction.

Consequently, $|u|_{X(a)}$ blows up as $a\to t_*-$. This shows (iv). \hfill $\square$
\section{Stability and instability results}
In this section we will assume that $f(0)=0$ and 
study the stability of the zero function for the semilinear problem (\ref{subproblem}).
We will further restrict ourselves to the case where the coefficient matrix $A$ does not depend
on time $t$. Recall that $\lambda_1>0$ denotes the first eigenvalue
of the negative Dirichlet-Laplacian $(-\Delta_D)$ in $L_2(\Omega)$.
If $A$ is also symmetric, by $\lambda_*>0$ we mean the smallest eigenvalue of the operator
$\mathcal{L}v=-\mbox{div}\big(A(x)\nabla v\big)$ (with Dirichlet boundary condition) in $L_2(\Omega)$.
We have the following stability result.
\begin{satz} \label{stabilitysubdiffusion}
Let $\Omega$ be a bounded domain in $\iR^N$ and $f\in C^{1-}(I)$, where $I$ is an open interval in $\iR$ containing $0$. Let $f(0)=0$ and assume that $f$ is differentiable at $0$. Suppose that the condition
($\mathcal{H}$) is satisfied and that $A$ is independent of $t$.
Let $l\in L_{1,\,loc}(\iR_+)$ be such that $(k,l)\in \mathcal{PC}$. Assume further that one of the two following {stability conditions}
is satisfied.
\begin{itemize}
\item [(a)] $A$ is also symmetric and $f'(0)<\lambda_*$.
\item [(b)] $f'(0)<\nu \lambda_1$.
\end{itemize}

Then $0$ is stable in the following sense: for any $\varepsilon>0$ there exists a $\delta>0$ such that
whenever $|u_0|_{L_\infty(\Omega)}\le \delta$ the problem (\ref{subproblem}) admits a global solution
$u$ satisfying
\begin{equation} \label{stabprop}
|u|_{L_\infty((0,\infty)\times \Omega)}\le \varepsilon.
\end{equation}
Moreover, if $|u_0|_{L_\infty(\Omega)}\le \delta$ we also have
\begin{equation} \label{relaxstab}
|u(t,x)|\le C\delta s_{\varepsilon_1}(t),\quad \mbox{a.a.}\;(t,x)\in (0,\infty)\times \Omega,
\end{equation}
for some $\varepsilon_1>0$ and some $C>0$ independent of $u$ and $u_0$.
In particular, if in addition $l\notin L_1(\iR_+)$,
then $0$ is even asymptotically stable, that is, $0$ is stable and
\[
|u|_{L_\infty((a,\infty)\times \Omega)}\to 0\quad \mbox{as}\;a\to \infty.
\]
\end{satz}
{\em Proof.} We give the argument in the case where the stability condition b) is satisfied.
The proof is the same for a); one only has to replace $\nu\lambda_1$ by $\lambda_*$ in the
subsequent formulas.

Given $\varepsilon>0$ we put $\varepsilon_0=(\nu \lambda_1-f'(0))/2$. Since $f$ is differentiable at $0$, there is $\rho>0$ such that $[-\rho,\rho]\subset I$ and
\begin{equation} \label{eps0ungl}
|f(y)-f'(0)y|=|f(y)-f(0)-f'(0)y|\le \varepsilon_0|y|\quad \mbox{for all}\,y\in [-\rho,\rho].
\end{equation}
Let $\delta$ be a number in the interval $(0,\rho/2]$ which will be fixed later.

Suppose that $u_0\in L_\infty(\Omega)$ with $|u_0|_\infty\le \delta$. Let $u$ be the corresponding solution
of problem (\ref{subproblem}) with maximal interval of existence $[0,t_*(u_0))$.
Let $\tau$ be the first exit time of $u$ for the interval $[-\rho,\rho]$, that is
\[
\tau:=\sup\{t_1\in (0,t_*(u_0)):\,|u(t,x)|\le \rho\,\;\;\mbox{a.a.}\;(t,x)\in (0,t_1)\times \Omega\}.
\]
From Theorem \ref{wellposedsubdiffusion} (ii) we know that $\tau>0$, since $\delta\le \rho/2$.

Next, suppose that $\tau<\infty$ and let $t_1\in (0,\tau)$. By Lemma \ref{pospartsub}, we have for the positive part of $u$ that
\begin{align*}
\partial_t\big(k\ast[u_+-(u_0)_+]\big)-\mbox{div}\big(A(x)\nabla [u_+]\big) & \le
f(u)\chi_{\{u\ge 0\}},\quad (t,x)\in \Omega_{t_1},
\nonumber\\
u_+ & = 0,\quad (t,x)\in \Gamma_{t_1}, \\
u_+|_{t=0} & = (u_0)_+,\quad x\in \Omega,\nonumber
\end{align*}
in the weak sense. Since $|u(t,x)|\le \rho$ for a.a.\ $(t,x)\in \Omega_{t_1}$, we
may use (\ref{eps0ungl}) and the relation $2\varepsilon_0=\nu\lambda_1-f'(0)$ to estimate
as follows (recall that $f(0)=0$).
\[
f(u)\chi_{\{u\ge 0\}}=f(u_+)\le f'(0)u_++\varepsilon_0 u_+= (\nu\lambda_1-\varepsilon_0)u_+.
\]
Setting $v_0=(u_0)_+$ we see that
\begin{align*}
\partial_t\big(k\ast[u_+-v_0]\big)-\mbox{div}\big(A(x)\nabla [u_+]\big) & \le
(\nu\lambda_1-\varepsilon_0)u_+,\quad (t,x)\in \Omega_{t_1},
\nonumber\\
u_+ & = 0,\quad (t,x)\in \Gamma_{t_1}, \\
u_+|_{t=0} & =  v_0,\quad x\in \Omega,\nonumber
\end{align*}
in the sense of a weak subsolution.

Now let $v$ be the bounded weak solution of 
\begin{align*}
\partial_t\big(k\ast[v-v_0]\big)-\mbox{div}\big(A(x)\nabla v\big) & =
(\nu\lambda_1-\varepsilon_0)v,\quad (t,x)\in \Omega_{t_1},
\nonumber\\
v & = 0,\quad (t,x)\in \Gamma_{t_1}, \\
v|_{t=0} & =  v_0,\quad x\in \Omega,\nonumber
\end{align*}
cf.\ Theorem \ref{linearestimate}. By the comparison principle, Theorem \ref{comp}, we have
$u_+\le v$ a.e.\ in $\Omega_{t_1}$. On the other hand, we know from Theorem \ref{linearestimate}
that there exists $\varepsilon_1>0$ and $\tilde{C}\ge 1$, both independent of $v$, $v_0$ and $t_1$, such that
\[
|v(t,x)|\le \tilde{C} s_{\varepsilon_1}(t)|v_0|_{L_\infty(\Omega)},\quad \mbox{a.a.}\,(t,x)\in\Omega_{t_1}.
\]
Thus
\[
|u_+(t,x)|\le \tilde{C} s_{\varepsilon_1}(t)|(u_0)_+|_{L_\infty(\Omega)},\quad \mbox{a.a.}\,(t,x)\in\Omega_{t_1}.
\]

Concerning the negative part of $u$, we set $\tilde{f}(y)=-f(-y)$ for $-y\in I$ and multiply the equation for $u$ by $-1$, thereby getting
\begin{align*}
\partial_t\big(k\ast[(-u)-(-u_0)]\big)-\mbox{div}\big(A(x)\nabla (-u)\big) & =\tilde{f}(-u),\quad (t,x)\in\Omega_{t_1},
\nonumber\\
-u & =0,\quad (t,x)\in\Gamma_{t_1},\\
(-u)|_{t=0} & = -u_0,\quad x\in \Omega.\nonumber
\end{align*}
We then proceed as above, now applying Lemma \ref{pospartsub} to $(-u)_+$. Note that
$\tilde{f}'(0)=f'(0)$ and thus by using (\ref{eps0ungl}) we have
\[
\tilde{f}(-u) \chi_{\{-u\ge 0\}}=\tilde{f}((-u)_+)\le f'(0)(-u)_++\varepsilon_0 (-u)_+= (\nu\lambda_1-\varepsilon_0)(-u)_+.
\]
By the same argument as above we now obtain
\[
|(-u)_+(t,x)|\le \tilde{C} s_{\varepsilon_1}(t)|(-u_0)_+|_{L_\infty(\Omega)},\quad \mbox{a.a.}\,(t,x)\in\Omega_{t_1}.
\]

Combining the estimates for the positive and negative part of $u$ yields
\begin{equation} \label{mainstab}
|u(t,x)|\le \tilde{C} s_{\varepsilon_1}(t)|u_0|_{L_\infty(\Omega)}\le \delta \tilde{C} s_{\varepsilon_1}(t)
,\quad \mbox{a.a.}\,(t,x)\in\Omega_{t_1}.
\end{equation}
Recall that $\tilde{C}\ge 1$. Choosing
\[
\delta=\tilde{C}^{-1}\,\min\big\{\varepsilon,\frac{\rho}{2}\big\}
\]
it follows from (\ref{mainstab}) that
\begin{equation} \label{ubound}
|u(t,x)|\le \min\big\{\varepsilon,\frac{\rho}{2}\big\},\quad \mbox{a.a.}\;(t,x)\in\Omega_{t_1}.
\end{equation}
Since $t_1<\tau$ was arbitrary, it follows that the estimate in (\ref{ubound}) even holds in $\Omega_{\tau}$.
By Theorem \ref{wellposedsubdiffusion} (iv) ($f|_{[-\rho,\rho]}$ can be extended to a function belonging to $C^{1-}(\iR)$) 
it is clear that $\tau<t_*(u_0)$. Knowing that $|u|\le \rho/2$ a.e.\ in $\Omega_{\tau}$ we can argue as in Step 3 in the
proof of Theorem \ref{wellposedsubdiffusion} to see that there exists $\tilde{\tau}\in (\tau,t_*(u_0))$ such that
$|u|\le \rho$ a.e.\ in $\Omega_{\tilde{\tau}}$. This contradicts the definition of $\tau$, so $\tau$ cannot be 
finite (as we assumed above). Theorem \ref{wellposedsubdiffusion} (iv) then implies $t_*(u_0)=\infty$.
Once we know this, (\ref{mainstab})
and  (\ref{ubound}) hold with $t_1$ being replaced by $\infty$ ($\Omega_\infty:=(0,\infty)\times \Omega$). 
In particular, (\ref{stabprop}) is satisfied.

Finally, if $l\notin
L_1(\iR_+)$ then we know from \cite[Lemma 6.1]{VZ} that $s_{\varepsilon_1}(t)\to 0$ as $t\to\infty$,
thereby proving the last assertion of the theorem. \hfill $\square$

${}$

\noindent We come now to an instability result. We will assume that $A$ does not depend
on time $t$ and is symmetric.
\begin{satz} \label{instabilitysubdiffusion}
Let $\Omega$ be a bounded domain in $\iR^N$ and $f\in C^{1-}(I)$, where $I$ is an open interval in $\iR$ containing $0$. Let $f(0)=0$ and assume that $f$ is differentiable at $0$. Suppose that the condition
($\mathcal{H}$) is satisfied and that $A$ is independent of $t$ and symmetric.
Assume that $(k,l)\in \mathcal{PC}$ with $l\notin L_1(\iR_+)$.
Suppose further that the {\em instability
condition}
\[
f'(0)>\lambda_*
\]
is fulfilled. Then $0$ is unstable.
\end{satz}
{\em Proof.} Fix $\varepsilon_0\in (0,f'(0)-\lambda_*)$. Since $f'(0)$ exists there is a $\rho>0$ such that
(\ref{eps0ungl}) is satisfied. Suppose that $0$ is stable. Then there exists $\delta>0$ such that for any $u_0\in L_\infty(\Omega)$ with $|u_0|_{L_\infty(\Omega)}\le \delta$ the corresponding
solution $u$ of (\ref{subproblem}) exists globally and $|u|_{L_\infty(0,\infty)\times \Omega)}\le \rho$. We choose $u_0\equiv \delta$. Appealing to Theorem  \ref{wellposedsubdiffusion} (iii), the solution $u$ is nonnegative.

Let $\psi\in \oH^1_2(\Omega)$ be the positive eigenfunction to the eigenvalue $\lambda_*$ with $|\psi|_{L_1(\Omega)}=1$, see e.g.\ \cite{GilTrud}. Fix $t_1>0$ and test
the time-regularized problem for $u$ with $\psi$. Using the eigenfunction property of $\psi$ this yields
\begin{align*}
\int_\Omega \Big(\partial_t \big[k_n\ast \big(u\psi-u_0\psi)\big]+\lambda_* u\psi\Big)\,dx=\int_\Omega f(u)\psi\,dx +\zeta_n(t),\quad \mbox{a.a.}\;t\in (0,t_1),
\end{align*}
where
\[
\zeta_n(t)=\int_\Omega \Big( \big(A\nabla u-h_n\ast [A\nabla u]|\nabla \psi\big)+\big(h_n\ast f(u)-f(u)\big)\psi \Big)\,dx.
\]

In view of (\ref{eps0ungl}), we have
\[
f(u)\ge f'(0)u-\varepsilon_0|u|=\big(f'(0)-\varepsilon_0)u.
\]
Setting $W(t)=\int_\Omega \psi u(t,x)dx$ and $W_0=\int_\Omega u_0 \psi\,dx=\delta$, it follows that
\begin{equation} \label{instabW}
\partial_t \big[k_n\ast (W-W_0)\big](t)\ge (f'(0)-\lambda_*-\varepsilon_0)W(t)+\zeta_n(t),\quad \mbox{a.a.}\;t\in (0,t_1).
\end{equation}

Convolving (\ref{instabW}) with $l$ and sending $n\to \infty$, the term involving
$\zeta_n$ drops ($|\zeta_n|_{L_1((0,t_1))}\to 0$), and after taking a subsequence if necessary we obtain
\[
W(t)\ge W_0+(f'(0)-\lambda_*-\varepsilon_0)(l\ast W)(t)\ge W_0=\delta>0,\quad \mbox{a.a.}\;t\in (0,t_1).
\]
Thus $W(t)$ is bounded away from zero. Returning to (\ref{instabW}), we
set $\kappa=f'(0)-\lambda_*-\varepsilon_0$,
divide
the inequality by $W(t)$ and apply Lemma \ref{convexFI} with $H(y)=-\log y$, $y>0$, to the
result
\begin{align*}
\partial_t \big[k_n\ast (\log W-\log W_0)\big](t)\ge \kappa+\frac{\zeta_n(t)}{W(t)},\quad \mbox{a.a.}\;t\in (0,t_1).
\end{align*}
Convolving next with $l$ and sending $n\to \infty$, the term involving $\zeta_n$
drops again and we get
\[
\log W(t)\ge \log W_0+\kappa(1\ast l)(t),\quad \mbox{a.a.}\;t\in (0,t_1).
\]
Since $t_1>0$ was arbitrary, this in turn implies
\begin{equation} \label{instabfin}
W(t)\ge W_0 e^{\kappa (1\ast l)(t)}=\delta e^{\kappa (1\ast l)(t)},\quad \mbox{a.a.}\;t>0.
\end{equation}
Since $l\notin L_1(\iR_+)$ and $\kappa>0$, the right-hand side of (\ref{instabfin}) becomes infinite as $t\to \infty$, which contradicts
\[
W(t)=\int_\Omega u(t,x)\psi(x)\,dx\le \rho,\quad \mbox{a.a.}\;t>0.
\]
The theorem is proved. \hfill $\square$
\section{Blowup}
\subsection{The purely time-dependent case}
We consider first the problem
\begin{equation} \label{nonlinprob}
\frac{d}{dt}\,\big(k\ast [u-u_0]\big)(t)=f(u(t)),\quad t>0,\;\;u(0)=u_0.
\end{equation}
Here $k$ is of type $\mathcal{PC}$. If we assume that $f:I\rightarrow \iR$ is locally Lipschitz continuous on the open interval $I\subset \iR$ and $u_0\in I$, then
(\ref{nonlinprob}) possesses a unique solution $u$ on a maximal interval of existence
$[0,t_*(u_0))$ with $u|_{[0,a]}\in L_\infty((0,a))$ for all $a<t_*(u_0)$. This can
be shown by similar arguments as in the Steps 1-3 in the proof of Theorem \ref{wellposedsubdiffusion}. It is not difficult to see, that the solution even has more regularity, it belongs to $H^1_1((0,a))$ for all $a<t_*(u_0)$.
Moreover, if $f(I)\subset [0,\infty)$ then $u(t)\ge u_0$ for all $t\in [0,t_*(u_0))$. Note that (\ref{nonlinprob}) is
equivalent to the Volterra equation
\[
u(t)=u_0+\big(l \ast  f(u)\big)(t),\quad t\ge 0.
\]
Such equations are studied thoroughly in the monograph \cite{GLS}.

The basic theorem on blowup of solutions to (\ref{nonlinprob}) is the following.
\begin{satz} \label{basicblowup}
Let $(k,l)\in \mathcal{PC}$, $\alpha_0,\alpha\in \iR$ with $\alpha_0<\alpha$. Let $f:(\alpha_0,\infty)\rightarrow \iR$
be locally Lipschitz continuous and $f|_{[\alpha,\infty)}:[\alpha,\infty)\rightarrow (0,\infty)$ be nondecreasing.
Assume further that
\begin{equation} \label{blowupcond}
\int_{\alpha}^\infty \frac{dr}{f(r)}<\infty.
\end{equation}
Then the following statements hold true.

(i) If $l\notin L_1(\iR_+)$ then for any $u_0\in [\alpha,\infty)$ the solution of (\ref{nonlinprob}) (in the class described above) blows up in finite time.

(ii) If $l\in L_1(\iR_+)$ then there exists $\beta \ge \alpha$ such that for any $u_0\in [\beta,\infty)$ the solution of (\ref{nonlinprob}) blows up in finite time.
\end{satz}
{\em Proof.} We will proceed by formal estimates, which can be made
rigorous by regularizing the problem in time. More precisely, convolving (\ref{nonlinprob}) with
$h_n$, $n\in \iN$ (cf.\ (\ref{hndef})), the kernel $k$ is replaced by $k_n$ (cf.\ (\ref{kndef})), which is admissible in
Lemma \ref{convexFI}. This is the same trick we used already in the previous sections.

Let $u_0\ge \alpha$ and suppose the solution exists globally. Multiplying the equation by $f(u(t))^{-1}$ then gives
\[
\frac{1}{f(u(t))}\,\frac{d}{dt}\,k\ast (u-u_0)=1,\quad t>0.
\]
Define
\[
F(y)=\int_{u_0}^y \frac{dr}{f(r)},\quad y\ge u_0.
\]
Then $F'(y)=f(y)^{-1}$ and $F$ is concave, since $f$ is nondecreasing. Furthermore, $F(u_0)=0$. By Lemma \ref{convexFI} it follows that
\[
\frac{d}{dt}\,\big(k\ast F(u)\big)\ge 1,\quad t>0.
\]
Convolving this inequality with the nonnegative kernel $l$ yields
\begin{equation} \label{lowerbound}
F(u(t))\ge (1\ast l)(t) ,\quad t\ge 0.
\end{equation}

In the case $l\notin L_1(\iR_+)$, the right-hand side of (\ref{lowerbound}) becomes infinite as $t\to \infty$. On the other
hand, the assumption (\ref{blowupcond}) implies that the left-hand side of (\ref{lowerbound}) stays bounded as $t\to \infty$, a contradiction. Hence $u$ does not exist globally, which means we have blowup in finite time for all $u_0\ge \alpha$.

In the case $l\in L_1(\iR_+)$ we choose $\beta\ge \alpha$ so large that
\[
\int_{\beta}^\infty \frac{dr}{f(r)}< |l|_{L_1(\iR_+)}.
\]  
Then for $u_0\ge \beta$ we deduce from (\ref{lowerbound}) that
\[
(1\ast l)(t)\le F(u(t))\le \int_{\beta}^\infty \frac{dr}{f(r)}< |l|_{L_1(\iR_+)}.
\]
Sending $t\to \infty$ leads to a contradiction.   \hfill $\square$
\begin{bemerk1} \label{blowupremark}
{\em
Inspection of the proof of Theorem \ref{basicblowup} shows that the statements of the previous theorem remain true for {\em weak supersolutions} of (\ref{nonlinprob}), that is, for $u:(0,t_*)\rightarrow \iR$ such that for all $a\in (0,t_*)$, we have $u|_{(0,a)}\in L_\infty((0,a))$ and
\begin{equation} \label{weaksuper1}
\int_0^a -\dot{\varphi}(t)\,\big(k\ast (u-u_0)\big)(t)\,dt\ge \int_0^a \varphi(t)f(u(t))\,dt,
\end{equation}
for all nonnegative $\varphi\in C^1([0,a])$ satisfying $\varphi(a)=0$. Note that
(\ref{weaksuper1}) is equivalent to
\[
\partial_t\big(k_n\ast (u-u_0)\big)(t)=\big(h_n\ast f(u)\big)(t),\quad\mbox{for a.a.}\;\,t\in (0,a),\,\mbox{and all}\;n\in \iN.
\]
Here $h_n$ is defined via (\ref{hndef}) as before.
}
\end{bemerk1}
\subsection{The PDE case}
We consider again the nonlocal PDE problem (\ref{subproblem}). We will assume that $A$ is independent of $t$ and
symmetric. Let again $\lambda_*>0$ denote the smallest eigenvalue of the operator $\mathcal{L}v=-\mbox{div}\big(A(x)\nabla v\big)$ (with Dirichlet boundary condition) and let
$\psi\in \oH^1_2(\Omega)$ denote the corresponding positive eigenfunction with $|\psi|_{L_1(\Omega)}=1$.
The proof of the following blowup result uses the eigenfunction method due to Kaplan \cite{Kap}, which is well known in the
classical parabolic case, see also \cite[Section 17]{QS}. 
\begin{satz} \label{blowupPDE}
Let $\Omega\subset \iR^N$ be a bounded domain and
$(k,l)\in \mathcal{PC}$. Suppose that the condition
($\mathcal{H}$) is satisfied and that $A$ is independent of $t$ and symmetric. Let $\alpha_0\in \iR$ and $f:(\alpha_0,\infty)\rightarrow \iR$ be a convex $C^1$-function. Suppose that
there exists $\alpha>\max\{0,\alpha_0\}$ such that $f(y)>0$ for all $y\ge \alpha$ and
\begin{equation} \label{supercond}
\int_\alpha^\infty \frac{dy}{f(y)}<\infty.
\end{equation}
Then there exists $M=M(\lambda_*,f,l)>0$ such that for any nonnegative $u_0\in L_\infty(\Omega)$
satisfying
\[
\int_\Omega u_0\psi\,dx \ge M
\]
the corresponding (weak) solution $u(t;u_0)$ of (\ref{subproblem}) blows up in finite time.
\end{satz}
{\em Proof.} We proceed again by formal estimates. Let $u_0\in L_\infty(\Omega)$ be nonnegative and $u$ be the corresponding solution of (\ref{subproblem}).
Let $T\in (0,t_*(u_0))$ be arbitrarily fixed. Taking $\psi$ as test-function
and setting $W(t)=\int_\Omega \psi u(t,x)dx$ and $W_0=\int_\Omega u_0 \psi\,dx$
we obtain (cf.\ the proof of Theorem \ref{instabilitysubdiffusion})
\begin{align*}
\partial_t \big(k\ast [W-W_0]\big)(t)+\lambda_*W(t)=\int_\Omega f(u)\psi \,dx,\quad
\mbox{a.a.}\;t\in (0,T).
\end{align*}
Since $f$ is convex and $|\psi|_{L_1(\Omega)}=1$, Young's inequality yields
\[
\int_\Omega f(u)\psi \,dx\ge f(W),
\]
and thus
\begin{align} \label{Wproblem}
\partial_t \big(k\ast [W-W_0]\big)(t)+\lambda_*W(t)\ge f(W(t)),\quad
\mbox{a.a.}\;t\in (0,T).
\end{align}
From the convexity of $f$ and assumption (\ref{supercond}) it follows that there exists $\alpha_1> \alpha$ such that $f(y)\ge 2\lambda_* y$ and $f'(y)\ge 0$ for all $y\ge \alpha_1$ (see also the proof of Theorem 17.3 in \cite{QS}).

Suppose now that $W_0\ge \alpha_1$. We claim that (\ref{Wproblem}) implies that
$W(t)\ge \alpha_1$ for a.a.\ $t\in (0,T)$. In fact, letting $M=|f(W)|_{L_\infty((0,T))}$ we can argue similarly as in the
proof of Theorem \ref{comp} to see that $W(t)\ge V(t)$ for a.a.\ $t\in (0,T)$ where $V\in C([0,T])$ solves the problem
\[
\partial_t \big(k\ast [V-W_0]\big)(t)+\lambda_*V(t)=-M,\;\;t\in (0,T),\quad V(0)=W_0.
\]
The solution $V$ is given by
\[
V(t)=s_{\lambda_*}(t)\,W_0-(1\ast r_{\lambda_*})(t)\,M,\quad t\in [0,T],
\]
cf.\ also the proof of Theorem \ref{linearestimate}. By continuity of $f$, there exists $\tilde{\alpha}_1\in [\alpha,\alpha_1)$
such that $f(y)\ge \lambda_* y$  for all $y\ge \tilde{\alpha}_1$. Recall that $V(0)=W_0\ge \alpha_1$.
Thus, by continuity of $V$, there exists $\delta\in (0,T]$ such that $V(t)\ge \tilde{\alpha}_1$
for all $t\in [0,\delta]$. This implies $W(t)\ge \tilde{\alpha}_1$
for a.a.\ $t\in (0,\delta)$, that is, $f(W(t))\ge \lambda_* W(t)$  for a.a.\ $t\in (0,\delta)$. Applying this estimate in (\ref{Wproblem}) and convolving the resulting inequality with the kernel $l$ yields $W(t)\ge W_0$ for a.a.\ $t\in (0,\delta)$.
Setting
\[
\delta_1:=\sup\{s\in (0,T):\,W(t)\ge W_0\,\;\;\mbox{for a.a.}\;t\in (0,s)\},
\]
we already know that $\delta_1>0$. Suppose that $\delta_1<T$. For $t\in (\delta_1,T)$ we may shift the time
as in Step 3 of the proof of Theorem \ref{wellposedsubdiffusion} by
setting $s=t-\delta_1$ and $\tilde{W}(s)=W(s+\delta_1)$, $s\in (0,T-\delta_1)$. By positivity of $W-W_0$ on $(0,\delta_1)$
and since $k$ is nonincreasing,
we have formally 
\begin{equation} \label{ts}
\partial_s\big(k\ast[\tilde{W}-W_0]\big)(s)\ge \partial_t\big(k\ast[W-W_0]\big)(s+\delta_1),\quad 
\mbox{a.a.}\;s\in  (0,T-\delta_1).
\end{equation}
This time-shifting property can be already found in \cite[Section 3.1]{ZwH} in the time fractional situation.
Note that the rigorous statement/argument uses the time-regularized version of the problem, where $k$ is replaced with the
more regular and nonincreasing kernel $k_n$. From (\ref{Wproblem}) and (\ref{ts}) we deduce that
\[
\partial_s \big(k\ast [\tilde{W}-W_0]\big)(s)+\lambda_*\tilde{W}(s)\ge f(\tilde{W}(s)),\quad
\mbox{a.a.}\;s\in (0,T-\delta_1),
\]
in the weak sense. So we may repeat the argument from above to see that there exists $\tilde{\delta}\in (0,T-\delta_1]$
such that $\tilde{W}(s)\ge W_0$ for a.a.\ $s\in (0,\tilde{\delta})$. This leads to a contradiction to the definition of $\delta_1$.
Hence, the assumption $\delta_1<T$ was not true. This proves the claim.

Knowing that $W(t)\ge W_0\ge \alpha_1$ for a.a.\ $t\in (0,T)$ it follows from (\ref{Wproblem}) that
\begin{equation} \label{Wsuper}
\partial_t \big(k\ast [W-W_0]\big)(t)\ge \frac{1}{2}\,f(W(t)),\quad
\mbox{a.a.}\;t\in (0,T).
\end{equation}
Since $f$ is nondecreasing on $[\alpha_1,\infty)$ and $\int_{\alpha_1}^\infty \frac{dy}{f(y)}<\infty$, we are in the situation of Remark \ref{blowupremark},
which says that there is some $M\ge \alpha_1$ depending only on $\alpha_1, l, f$ such that for any nonnegative $u_0\in L_\infty(\Omega)$
satisfying $W_0\ge M$, the function $W(t)$ satisfying (\ref{Wsuper}) blows up in finite time, and thus the same holds for the
(weak) solution $u(t;u_0)$ of (\ref{subproblem}).  \hfill $\square$


$\mbox{}$
{\footnotesize

\noindent {\bf Vicente Vergara}, Departamento de Matem\'{a}ticas, Universidad de La Serena, Avenida
Cisternas 1200, La Serena, Chile, and
Universidad de Tarapac\'{a}, Avenida General Vel\'{a}squez 1775, Arica, Chile, E-mail: vvergaraa@uta.cl

$\mbox{}$

\noindent {\bf Rico Zacher}, Institute of Applied Analysis, University of Ulm, 89069 Ulm, Germany, E-mail: rico.zacher
@uni-ulm.de

}

\end{document}